\documentclass{amsart}
\usepackage{amssymb}
\usepackage{amsfonts}
\usepackage{latexsym}
\usepackage[all]{xy}
\usepackage{amscd}

\newtheorem{theorem}{Theorem}[section]
\newtheorem{lemma}[theorem]{Lemma}
\newtheorem{proposition}[theorem]{Proposition}
\newtheorem{corollary}[theorem]{Corollary}
\theoremstyle{definition}
\newtheorem{definition}[theorem]{Definition}
\newtheorem{example}[theorem]{Example}

\newtheorem{remark}[theorem]{Remark}

\newcommand{\id}{\text{id}}

\newcommand{\lieg}{\mathfrak{g}}
\newcommand{\te}{\theta}

\newcommand{\cA}{\mathcal{A}}
\newcommand{\cP}{\mathcal{P}}
\newcommand{\Wtil}{\tilde{W}}

\newcommand{\kruisje}[1]{\, \mbox{$_{#1}$}\hspace{-.2ex}\mbox{$\ltimes$}
\,}

\newcommand{\alh}{\hat{\alpha}}
\newcommand{\recht}{\to}
\newcommand{\io}{\id}
\newcommand{\cU}{\mathcal{U}}
\newcommand{\cV}{\mathcal{V}}
\newcommand{\al}{\alpha}
\newcommand{\be}{\beta}
\newcommand{\ot}{\otimes}
\newcommand{\om}{\omega}
\newcommand{\de}{\Delta}
\newcommand{\sde}{\delta}
\newcommand{\Mh}{\hat M}
\newcommand{\deh}{\hat\Delta}
\newcommand{\deop}{\Delta^{op}}

\newcommand{\C}{\mathcal{C}}

\newcommand{\cO}{\mathcal{O}}

\begin{document}
\title{The bicrossed product construction for locally compact
quantum groups}
\author{Leonid Vainerman}
\address{Dept. of Math., University of Caen, B.P. 5186, 14032 Caen
Cedex, France; E.mail: vainerma@math.unicaen.fr}
\date{October 3, 2005}
\begin{abstract}
The cocycle bicrossed product construction allows certain freedom
in producing examples of locally compact quantum groups. We give
an overview of some recent examples of this kind having remarkable
properties.
\end{abstract}
\maketitle

\begin{section}
{Introduction}

The initial motivation to introduce objects which are more general
than locally compact (l.c.) groups was to extend classical
harmonic analysis on these groups, including the Fourier transform
and the Pontrjagin duality. Given an abelian l.c. group $G$, the
set $\hat G$ of its unitary continuous characters is again an
abelian l.c. group - the dual group of $G$. The Fourier transform
maps functions on $G$ to functions on $\hat G$, and the Pontrjagin
duality theorem claims that the dual of $\hat G$ is isomorphic to
$G$. If $G$ is not abelian, the set of its characters is too
small, and one should use instead the set $\hat G$ of (classes of)
its unitary irreducible representations and their matrix
coefficients. For compact groups, this leads to the Peter-Weyl
theory; the corresponding duality is due to T. Tannaka and M.G.
Krein. In this case, $\hat G$ is not a group, but carry a
structure of a {\it block-algebra} or a {\it Krein algebra}
\cite{HR}; however, this structure allows to reconstruct the
initial group. Such a non-symmetric duality was later established
by W.F. Stinespring for unimodular groups, and by P. Eymard and T.
Tatsuuma for general l.c. groups.

In order to restore the symmetry of the duality, G.I. Kac \cite{K}
introduced in 1961 a category of {\it ring groups}. A ring group
is a Hopf-von Neumann algebra $(M,\de,S)$, with the
comultiplication $\Delta : M\to M\otimes M$ and the involutive
antipode $S : M\to M,\ S^2=\id$ equipped with a faithful normal
{\it trace} $\varphi$ compatible with $\Delta$ and $S$ and playing
the role of a Haar measure. If $M$ is commutative (resp., $\de$ is
co-commutative, i.e., $\sigma\circ\Delta=\Delta$, where
$\sigma:a\otimes b\mapsto b\otimes a$ is the usual flip in
$M\otimes M$), this ring group can be identified with the algebra
$L^\infty(G)$ (resp., group von Neumann algebra $\mathcal{L}(G)$),
where $G$ is a unimodular group. Thus, unimodular groups and their
duals are embedded into this category, and the duality constructed
by Kac extended those of Pontrjagin, Tannaka-Krein and
Stinespring.

The theory was completed in early 70-s by G.I. Kac and the author,
and independently by M. Enock and J.-M. Schwartz, in order to
cover all locally compact groups. They allowed $\varphi$ and
$\varphi\circ S$ to be different {\it weights} on $M$ playing
respectively the role of left and right Haar measures (for ring
groups $\varphi=\varphi\circ S$ was a trace), gave appropriate
axioms and extended the construction of the dual. These more
general objects are called {\it Kac algebras} \cite{ES}. L.c.
groups and their duals can be viewed respectively as commutative
and co-commutative Kac algebras, the corresponding duality covered
all versions of duality for such groups.

Concrete examples of ring groups, which were neither ordinary
groups nor their duals, were given in \cite{KP1} - \cite{KP3}.
According to V.G. Drinfeld \cite{D}, these were the first concrete
examples of what is now called quantum groups.

Quantum groups discovered by V.G. Drinfeld \cite{D} and others
gave new important examples of Hopf algebras, obtained by
deformation of universal enveloping algebras and of function
algebras on Lie groups. Their operator algebra versions did not
fit into the Kac algebra theory, because their antipodes were
neither involutive, nor bounded maps. This motivated strong
efforts to construct a more general theory, which would be as
elegant as that of Kac algebras but would also cover these new
examples. Important steps in this direction were made by: S.L.
Woronowicz \cite{W1} - \cite{W8} with his theory of compact
quantum groups and a series of important concrete examples of
compact and non-compact quantum groups; S. Baaj and G. Skandalis
\cite{B} - \cite{BS1} with their fundamental concept of a
multiplicative unitary; T. Masuda, Y. Nakagami and S.L. Woronowicz
\cite{MN}, \cite{MNW} who gave a set of axioms of so called
Woronowicz algebras; A. Van Daele \cite{VD2}, \cite{VD3} who
introduced an important notion of a multiplier Hopf algebra.
Finally, the theory of l.c. quantum groups was proposed by J.
Kustermans and S. Vaes \cite{KV1} - \cite{KV3}.

A (von Neumann algebraic) l.c. quantum group is a collection
$(M,\Delta, \varphi,\psi)$, where $M$ is a von Neumann algebra
equipped with a co-associative comultiplication $\Delta:M\to
M\otimes M$ and two {\it normal semi-finite faithful} (n.s.f.)
weights $\varphi$ and $\psi$ - right and left Haar measures. The
antipode is not explicitly present in this definition, but can be
constructed from the above data. Kac algebras, compact and
discrete quantum groups are special cases of a l.c. quantum group,
and all important concrete examples of operator algebraic quantum
groups fit into this framework. There is an equivalent
$C^*$-algebraic version of a l.c. quantum group.

A number of "isolated" examples of l.c. quantum groups can be
formulated in terms of generators of certain Hopf $*$-algebras and
commutation relations between them. It is much harder to represent
these generators as (typically, unbounded) operators acting on a
Hilbert space, to give a meaning to the relations of commutation
between them, to associate an operator algebra with the above
system of operators and commutation relations and to construct
comultiplication, antipode and invariant weights as applications
related to this algebra. There is no general approach to these
highly nontrivial problems, and one must design specific methods
in each specific case (see, for example, \cite{B}, \cite{KK},
\cite{KK1}, \cite{PW}, \cite{VW}, \cite{W1} - \cite{W8}).

A systematic approach to the construction of non-trivial Kac
algebras has been proposed in \cite{K1}. Given two finite groups,
$G_1$ and $G_2$, viewed as a co-commutative and a commutative Kac
algebra, $(\mathcal{L}(G_1), \Delta_1)$ and $(L^\infty(G_2),
\Delta_2)$ respectively, let us try to find a Kac algebra
$(M,\Delta)$ which makes the sequence
\begin{equation}\label{exact}
(L^\infty(G_2), \Delta_2)\to (M,\Delta)\to (\mathcal{L}(G_1),
\Delta_1)
\end{equation}
exact. Kac explained that: 1) $(M,\Delta)$ exists if and only if
$G_1$ and $G_2$ are subgroups of a group $G$ such that $G_1\cap
G_2= \{e\}$ and $G=G_1G_2$. Equivalently, $G_1$ and $G_2$ must act
on each other (as on sets), and these actions must be compatible
in certain sense. Remark, that later on, such pairs of groups were
also introduced by G.W. Mackey \cite{Mac} and by M. Takeuchi
\cite{Tak}. Following \cite{Tak}, we will say that $G_1$ and $G_2$
form a {\it matched pair}. 2) To get all possible $(M,\Delta)$
(they are called {\it extensions} of $(L^\infty(G_2), \Delta_2)$
by $(\mathcal{L}(G_1), \Delta_1)$), one must find all possible
2-cocycles for the above mentioned actions, compatible in certain
sense. Then \cite{K1} gives explicit construction of $(M,\Delta)$
(the cocycle bicrossed product construction). The famous
Kac-Paljutkin examples \cite{KP1} - \cite{KP3} are of this type.

Later on, both algebraic and analytic aspects of this construction
were studied by S. Majid \cite{Maj1} - \cite{Majbook} who gave
also a number of examples of operator algebraic quantum groups.
The bicrossed product construction for multiplicative unitaries
was done in \cite{BS}. A general theory of extensions of the form
(\ref{exact}), with l.c. $G_1$ and $G_2$, has been developed in
\cite{VV1} (in fact, \cite{VV1} treated general case, with $G_1$
and $G_2$ l.c. quantum groups, which we don't consider here).

In the recent years it became clear that this construction gives
certain freedom in producing examples of l.c. compact quantum
groups, see \cite{BSV} - \cite{BV2}, \cite{CKP1} - \cite{DQV},
\cite{F}, \cite{SV2} - \cite{VV2}. This allows to construct such
examples with prescribed special features. Some of these examples
have quite remarkable and unexpected properties, and we discuss
them briefly here. We recall in Section 2 the basic definitions
and results on the von Neumann algebraic version of the l.c.
quantum group theory and the cocycle bicrossed product
construction. In Section 3 we discuss the problem of regularity
for a multiplicative unitary related to a l.c. quantum group and
give, following \cite{BSV}, a surprising example of a
non-semi-regular l.c. quantum group. Section 4 is devoted to the
examples of l.c. quantum groups whose von Neumann algebras are
factors. Finally, we briefly discuss in Section 5 amenability and
Kac exact sequence for l.c. quantum groups coming from the cocycle
bicrossed product construction.

{\bf Acknowledgment.} The author expresses his gratitude to the
organizers of the International Conference on Harmonic Analysis
and Quantum Groups in Cochin, India (January, 2005) for their kind
invitation to present a talk there. He is grateful to S. Vaes for
many useful discussions on the bicrossed product construction and
on related topics, and also to Zh.-J. Ruan and A. Van Daele for
helpful remarks.
\end{section}
\begin{section}
{Preliminaries}
\begin{subsection}{Notations.}
Let us denote by $B(H)$ the algebra of all bounded linear
operators on a Hilbert space $H$, by $\ot$ the tensor product of
Hilbert spaces or von Neumann algebras and by $\Sigma$ (resp.,
$\sigma$) the flip map on it. If $H, K$ and $L$ are Hilbert spaces
and $X \in B(H \ot L)$ (resp., $X \in B(H \ot K), X \in B(K \ot
L)$), we denote by $X_{13}$ (resp., $X_{12},\ X_{23}$) the
operator $(1 \ot \Sigma^*)(X \ot 1)(1 \ot \Sigma)$ (resp., $X\ot
1,\ 1\ot X$) defined on $H \ot K \ot L$. When $H = H_1 \ot H_2$
itself is a tensor product of two Hilbert spaces, we switch from
the above leg-numbering notation with respect to $H \ot K \ot L$
to the one with respect to the finer tensor product $H_1 \ot
H_2\ot K \ot L$, for example, from $X_{13}$ to $X_{124}$. There is
no confusion here, because the number of legs changes.

Given a comultiplication $\de$, denote by $\deop$ the opposite
comultiplication $\sigma \de$. Our reference to the modular theory
of weights on von Neumann algebras is \cite{Stra}. Given a {\it
normal semi-finite faithful} (n.s.f.) weight $\theta$ on a von
Neumann algebra $N$, we denote:
$$
{\mathcal M}^+_\theta = \{ x \in N^+ \mid \theta(x) < + \infty \},
\ {\mathcal N}_\theta = \{ x \in N \mid x^*x \in {\mathcal
M}^+_\theta \}\ \text{and}\ {\mathcal M}_\theta =
\operatorname{span} {\mathcal M}^+_\theta \; .
$$
All l.c. groups considered in this paper supposed to be second
countable.
\end{subsection}
\subsection{L.c.\ quantum groups \cite{KV1}-\cite{KV3}}
A pair $(M,\de)$ is called a (von Neumann algebraic) l.c.\ quantum
group  when
\begin{itemize}
\item $M$ is a von Neumann algebra and $\de : M \to M \ot M$ is a
normal and unital $*$-homomorphism which is coassociative: $(\de
\ot \id)\de = (\id \ot \de)\de$. \item There exist n.s.f. weights
$\varphi$ and $\psi$ on $M$ such that
\begin{itemize}
\item $\varphi$ is left invariant in the sense that $\varphi \bigl( (\om \ot
\id)\de(x) \bigr) = \varphi(x) \om(1)$ for all $x \in {\mathcal
M}_{\varphi}^+$ and $\om \in M_*^+$,
\item $\psi$ is right invariant in the sense that $\psi \bigl( (\id \ot
\om)\de(x) \bigr) = \psi(x) \om(1)$ for all $x \in {\mathcal
M}_{\psi}^+$ and $\om \in M_*^+$.
\end{itemize}
\end{itemize}
Left and  right invariant weights are unique up to a positive
scalar.

Represent $M$ on the GNS Hilbert space of $\varphi$ and define a
unitary $W$ on $H \ot H$ by
$$
W^* (\Lambda(a) \ot \Lambda(b)) = (\Lambda \ot \Lambda)(\de(b)(a
\ot 1)) \quad\text{for all}\; a,b \in N_{\phi}\; .
$$
Here, $\Lambda$ denotes the canonical GNS-map for $\varphi$,
$\Lambda \ot \Lambda$ the similar map for $\varphi \ot \varphi$.
One proves that $W$ satisfies the {\it pentagonal equation}:
$W_{12} W_{13} W_{23} = W_{23} W_{12}$, and we say that $W$ is a
{\it multiplicative unitary}. The von Neumann algebra $M$ and the
comultiplication on it can be given in terms of $W$ respectively
as
$$M = \{ (\id \ot \om)(W) \mid \om \in B(H)_* \}^{-\sigma-strong*} \; $$
and $\de(x) = W^* (1 \ot x) W$, for all $x \in M$. Next, the l.c.\
quantum group $(M,\de)$ has an antipode $S$, which is the unique
$\sigma$-strongly* closed linear map from $M$ to $M$ satisfying
$(\id \ot \om)(W) \in {\mathcal D}(S)$ for all $\om \in B(H)_*$
and $S(\id \ot \om)(W) = (\id \ot \om)(W^*)$ and such that the
elements $(\id \ot \om)(W)$ form a $\sigma$-strong* core for $S$.
$S$ has a polar decomposition $S = R \tau_{-i/2}$, where $R$ is an
anti-automorphism of $M$ and $(\tau_t)$ is a strongly continuous
one-parameter group of automorphisms of $M$. We call $R$ the
unitary antipode and $(\tau_t)$ the scaling group of $(M,\de)$. We
have $\sigma (R \ot R) \de = \de R$, so $\varphi R$ is a right
invariant weight on $(M,\de)$ and we take $\psi:= \varphi R$.

Let us denote by $(\sigma_t)$ the modular automorphism group of
$\varphi$. There exist a number $\nu > 0$, called the scaling
constant, such that $\psi \, \sigma_t = \nu^{-t} \, \psi$ for all
$t \in \mathbb{R}$. Hence, we get the existence of a unique
positive, self-adjoint operator $\sde_M$ affiliated to $M$, such
that $\sigma_t(\sde_M) = \nu^t \, \sde_M$ for all $t \in
\mathbb{R}$ and $\psi = \varphi_{\sde_M}$. The operator $\sde_M$
is called the modular element of $(M,\de)$. If $\sde_M=1$ we call
$(M,\de)$ unimodular. The scaling constant can be characterized as
well by the relative invariance $\varphi \, \tau_t = \nu^{-t} \,
\varphi$.

For the dual l.c.\ quantum group $(\Mh,\deh)$ we have :
$$\Mh = \{(\om \ot \id)(W) \mid \om \in B(H)_*
\}^{-\sigma-strong*}$$ and $\deh(x) = \Sigma W (x \ot 1) W^*
\Sigma$ for all $x \in \Mh$. If we turn the predual $M_*$ into a
Banach algebra with product $\om \, \mu = (\om \ot \mu)\de$ and
define
$$\lambda: M_* \to \Mh : \lambda(\om) = (\om \ot \id)(W),$$ then $\lambda$ is a
homomorphism and $\lambda(M_*)$ is a $\sigma$-strongly* dense
subalgebra of $\Mh$. A left invariant n.s.f. weight $\hat \varphi$
on $\Mh$ can be constructed explicitly and the associated
multiplicative unitary $\hat{W} = \Sigma W^* \Sigma$.

Since $(\Mh,\deh)$ is again a l.c.\ quantum group, we can
introduce the antipode $\hat{S}$, the unitary antipode $\hat{R}$
and the scaling group $(\hat{\tau}_t)$ exactly as we did it for
$(M,\de)$. Also, we can again construct the dual of $(\Mh,\deh)$,
starting from the left invariant weight $\hat\varphi$. The bidual
l.c.\ quantum group $(\hat\Mh,\hat\deh)$ is isomorphic to
$(M,\de)$.

We denote by $(\hat\sigma_t)$ the modular automorphism group of
the weight $\hat\varphi$. The modular conjugations of the weights
$\varphi$ and $\hat\varphi$ will be denoted by $J$ and $\hat J$
respectively. Then it is worthwhile to mention that
$$
R(x) =\hat J x^* \hat J, \quad\text{for all} \; x \in M,
\qquad\text{and}\qquad \hat R(y) = J y^* J, \quad\text{for all}\;
y \in \Mh \; .$$ Let us mention important special cases of l.c.\
quantum groups.

a) {\it Kac algebras} \cite{ES}. $(M,\de)$ is a Kac algebra if and
only if $(\tau_t)=\id$ and $\sigma_t \, R = R \, \sigma_{-t}$ for
all $t \in \mathbb{R}$ \cite{ES}. Let $(\sigma'_t)$ be the modular
automorphism group of $\psi$. Since $\psi = \varphi R$, we get
$\sigma'_t \, R = R \, \sigma_{-t}$ for all $t \in \mathbb{R}$.
Hence $(M,\de)$ is a Kac algebra if and only if $(\tau_t)=\id$ and
$\sigma'=\sigma$ or if and only if $\sde_M$ is affiliated to the
center of $M$.

In particular, $(M,\de)$ is a Kac algebra if $M$ is commutative.
Then $(M,\de)$ is generated by a usual l.c.\ group $G :
M=L^{\infty}(G), (\de f)(g,h) = f(gh),$ $(Sf)(g) = f(g^{-1}),\
\varphi(f)=\int f(g)\; dg$, where $f\in L^{\infty}(G),\ g,h\in G$
and we integrate with respect to the left Haar measure $dg$ on
$G$. The right invariant weight $\psi$ is given by $\psi(f) = \int
f(g^{-1}) \; dg$. The modular element $\sde_M$ is given by the
strictly positive function $g \mapsto \sde_G(g)^{-1}$.

The von Neumann algebra $L^\infty(G)$ acts on $H=L^2(G)$ by
multiplication and $$(W_G\xi)(g,h)=\xi(g,g^{-1}h)$$ for all
$\xi\in H\ot H=L^2(G\times G)$. Then $\Mh=\mathcal L(G)$ is the
group von Neumann algebra generated by the operators
$(\lambda_g)_{g\in G}$ of the left regular representation of $G$
and $\deh(\lambda_g)=\lambda_g\ot\lambda_g$. Clearly,
$\deh^{op}:=\sigma\deh=\deh$, so $\deh$ is cocommutative.

b) A l.c.\ quantum group is called {\it compact} if
$\varphi(1)<+\infty$. A l.c.\ quantum group $(M,\de)$ is called
{\it discrete} if $(\Mh,\deh)$ is compact.
\subsection{Cocycle crossed and bicrossed products} An {\it action} of a
l.c.\ quantum group $(M,\de)$ on a von Neumann algebra $N$ is a
normal, injective and unital $*$-homomorphism $\al: N \to M \ot N$
such that $(\id \ot \al)\al(x) =  (\de \ot \id)\al(x)$ for all $x
\in N$. This generalizes the definition of an action of a l.c.\
group $G$ on a ($\sigma$-finite) von Neumann algebra $N$, as a
continuous map $G \to\operatorname{Aut} N : s \mapsto \al_s$ such
that $\al_{st}=\al_s \al_t$ for all $s,t\in G$. Indeed, putting
$M=L^\infty(G)$, one can identify $M \ot N$ with $L^\infty(G,N)$
and $M \ot M \ot N$ with $L^\infty(G \times G,N)$ and define the
above homomorphism $\al$ by $(\al(x))(s) = \al_{s^{-1}}(x)$. The
fixed point algebra of an action $\al$ is defined by
$N^\alpha=\{x\in N\mid \al(x)=1\ot x\}$.

A {\it cocycle} for an action of a l.c.\ group $G$ on a
commutative von Neumann algebra $N$ is a Borel map $u:G \times G
\to N$ such that $\al_r(u(s,t)) \; u(r,st) = u(r,s) \; u(rs,t)$
nearly everywhere. Then, putting $M=L^\infty(G)$, one can define a
unitary $\cU \in M \ot M \ot N$ by $\cU(s,t)=u(t^{-1},s^{-1})$
satisfying $$(\id \ot \id \ot \al)(\cU) (\de \ot \id \ot \id)(\cU)
= (1 \ot \cU) (\id \ot \de \ot \id \cU) \; .$$ The general
definition of a cocycle action of a l.c.\ quantum group on a von
Neumann algebra can be found in \cite{VV1}. The {\it cocycle
crossed product} $G \kruisje{\al,\cU} N$ is the von Neumann
subalgebra of $B(L^2(G)) \ot N$ generated by $$\al(N)\ \text{and}
\ \{ (\om \ot \id \ot \id)(\Wtil) \mid \om \in L^1(G) \} \; ,$$
where $\Wtil = (W_G \ot 1)\cU^*$. There exists a unique action
$\alh$ of $(\mathcal{L}(G),\deh)$ on $G \kruisje{\al,\cU} N$ such
that
$$
\alh(\al(x)) = 1 \ot \al(x)\ \text{for all}\ x \in N, \\
(\id \ot \alh)(\Wtil) = W_{G,12} \Wtil_{134} \; ,
$$
and for any n.s.f. weight $\gamma$ on $N$, we can define the {\it
dual} n.s.f. weight $\tilde\gamma$ on $G \kruisje{\al,\cU} N$:
$$\tilde\gamma = \gamma \al^{-1} \; (\hat\varphi \ot \id \ot \id) \alh.$$
\begin{definition}\label{mpg} (see \cite{BSV}) Let $G,G_1$ and $G_2$
be l.c.\ groups, and let a homomorphism $i: G_1 \recht G$ and an
anti-homomorphism $j:G_2 \recht G$ have closed images and be
homeomorphisms onto these images. Suppose that $i(G_1) \cap j(G_2)
= \{e\}$ and that the complement of $i(G_1)j(G_2)$ in $G$ has
measure zero. Then we call $(G_1,G_2\subset G)$ a matched pair of
l.c.\ groups.
\end{definition}
\begin{remark}\label{double}
The above mentioned group $G$ is called a {\it double crossed
product} of $G_1$ and $G_2$. The definition of a matched pair of
general l.c. quantum groups was given in \cite{VV1}, and the
double crossed product construction in this general case was
studied in \cite{BV2}. It is worthwhile to mention that this
construction contains the famous Drinfeld double construction for
quantum groups \cite{D}.
\end{remark}
The map $\te:G_1\times G_2\to G: (g,s)\mapsto i(g)j(s)$ is
automatically a Borel isomorphism, i.e., it induces an isomorphism
between $L^\infty(G_1 \times G_2)$ and $L^\infty(G)$ \cite{BSV}.
Hence, this data allows to construct as follows two actions: $\al$
of $G_1$ on $M_2=L^\infty(G_2)$ and $\be$ of $G_2^{op}$ on
$M_1=L^\infty(G_1)$, verifying certain compatibility relations.

Let $\Omega$ be the image of $\te$ and define the Borel
isomorphism
$$
\rho : G_1 \times G_2 \recht \Omega^{-1} : (g,s) \mapsto j(s)i(g)
\; .
$$
So $\cO=\theta^{-1}(\Omega \cap \Omega^{-1})$ and
$\cO'=\rho^{-1}(\Omega \cap \Omega^{-1})$ are Borel subsets of
$G_1 \times G_2$, with complements of measure zero, and $\rho^{-1}
\theta$ is a Borel isomorphism of $\cO$ onto $\cO'$. For all
$(g,s) \in \cO$ define $\be_s(g) \in G_1$ and $\al_g(s) \in G_2$
such that
$$
\rho^{-1}(\theta(g,s)) = (\be_s(g),\al_g(s)) \; .
$$
Hence we get $j \bigl( \al_g(s) \bigr) \; i\bigl( \be_s(g) \bigr)
= i(g)j(s)$ for all $(g,s) \in \cO$.
\begin{lemma}\label{42} (\cite{VV1}, Lemma 4.8)
Let $(g,s) \in \cO$ and $h \in G_1$. Then $(hg,s) \in \cO$ if and
only if $(h,\al_g(s)) \in \cO$, and in that case
$$\al_{hg}(s) = \al_h \bigl( \al_g(s) \bigr) \qquad\text{and}\qquad \be_s(hg) =
\be_{\al_g(s)}(h) \; \be_s(g) \; .$$ Let $(g,s) \in \cO$ and $t
\in G_2$. Then $(g,ts) \in \cO$ if and only if $(\be_s(g),t) \in
\cO$ and in that case
$$\be_{ts}(g) = \be_t \bigl( \be_s(g) \bigr) \qquad\text{and}\qquad \al_g(ts) =
\al_{\be_s(g)}(t) \; \al_g(s) \; .$$ Finally, for all $g \in G_1$
and $s \in G_2$ we have $(g,e) \in \cO$, $(e,s) \in \cO$, and
$$\al_g(e) = e, \quad \al_e(s) = s, \quad \be_s(e)=e\ \text{and}\ \be_e(g)=g \; .$$
\end{lemma}
This can be viewed as a definition of a matched pair of l.c.\
groups in terms of mutual actions.

The cocycles for the above actions can be introduced as measurable
maps $\cU  : G_1 \times G_1 \times G_2\to \mathbb{T}$ and $\cV  :
G_1 \times G_2 \times G_2\to \mathbb{T}$, where $\mathbb{T}$ is
the unit circle in $\mathbb{C}$, satisfying
$$
\cU(g,h,\al_k(s)) \; \cU(gh, k, s) = \cU(h,k,s) \; \cU(g, hk, s),
$$
\begin{equation}\label{cocU}
\cV(\be_s(g),t,r) \; \cV(g,s, rt) = \cV(g,s,t) \; \cV(g, ts, r) ,
\end{equation}
$$\cV(gh, s, t) \; \bar{\cU}(g,h,ts)  =$$
$$=\bar{\cU}(g,h,s) \;
\bar{\cU}(\be_{\al_h(s)}(g), \be_s(h) , t)\cV(g, \al_h(s),
\al_{\be_s(h)}(t)) \; \cV(h,s,t) \notag
$$
nearly everywhere. This gives a definition of a cocycle matched
pair of l.c.\ groups.

Fixing a cocycle matched pair of l.c.\ groups $G_1$ and $G_2$,
denoting $H_i=L^2(G_i)$ $(i=1,2)$, $H=H_1\ot H_2$ and identifying
$\cU$ and $\cV$ with unitaries in $M_1\ot M_1\ot M_2$ and in
$M_1\ot M_2\ot M_2$ respectively, define unitaries $W$ and $\hat
W$ on $H \ot H$ by
$$
\hat W = (\be \ot \io \ot \io) \bigl( (W_{G_1} \ot 1) \cU^* \bigr)
\; (\io \ot \io \ot \al) \bigl( \cV (1 \ot \hat W_{G_2})
\bigr)\;\text{and}\; W = \Sigma \hat W^* \Sigma \; .
$$
On the von Neumann algebra $M = G_1 \kruisje{\al,\cU}
L^{\infty}(G_2)$, let us define a faithful $*$-homomorphism
$$\de : M \recht B(H \ot H) : \de(z) = W^*(1 \ot z) W \ (\forall z\in M)$$
and denote by $\varphi$ the dual weight of the canonical left
invariant trace $\varphi_2$ on $L^{\infty}(G_2)$. Then, Theorem
2.13 of \cite{VV1} shows that $(M,\de)$ is a l.c.\ quantum group
with $\varphi$ as a left invariant weight, which we call the {\em
cocycle bicrossed product} of $G_1$ and $G_2$. One can also show
that its scaling constant is 1. The dual l.c.\ quantum group is
$(\Mh,\deh)$, where $\Mh= G_2 \kruisje{\al,\cU} L^{\infty}(G_1)$\
and $\deh(z) = \hat W^* (1 \ot z) \hat W$ for all $z \in \Mh$.

One can get explicit formulas for the modular operators and
conjugations of the left invariant weights, unitary antipodes,
scaling groups and modular elements of both $(M,\de)$ and its dual
in terms of $\alpha_g,\ \beta_s$, the cocycles and the modular
functions $\delta,\delta_1$ and $\delta_2$ of the l.c.\ groups
$G,G_1$ and $G_2$. These formulas imply
\begin{proposition} \label{charKac}
The l.c.\ quantum group $(M,\de)$ is a Kac algebra if and only if
$$
\sde\bigl( i(g \be_s(g)^{-1}) \bigr) \; \sde_1\bigl(g^{-1}
\be_s(g) \bigr) \; \sde_2\bigl( \al_g(s) s^{-1} \bigr) \ = 1\
\text{and}\ \frac{\sde_1 \bigl( \be_s(g) \bigr)}{\sde_1(g)} =
\frac{\sde_2 \bigl( \al_g(s) \bigr)}{\sde_2(s)} \; .
$$
\end{proposition}
\begin{corollary}\label{410}
If $\al$ or $\be$ is trivial, $(M,\de)$ and $(\Mh,\deh)$ are Kac
algebras.
\end{corollary}
\begin{corollary}\label{411}
If both $\al$ and $\be$ preserve modular functions and Haar
measures, then $(M,\de)$ and $(\Mh,\deh)$ are Kac algebras.
\end{corollary}
Remark that the conditions of this corollary are fulfilled if both
groups are discrete since any discrete group is unimodular and its
Haar measure is constant.
\begin{corollary}\label{hom}
If $(G_1,G_2\subset G)$ is a fixed matched pair of l.c.\ groups
and cocycles $\cU$ and $\cV$ satisfy (\ref{cocU}), we get a
cocycle bicrossed product $(M,\de)$. If one of these cocycle
bicrossed products is a Kac algebra, then all of them are Kac
algebras.
\end{corollary}
\begin{proof}
The necessary and sufficient conditions for $(M,\de)$ to be a Kac
algebra in Proposition~\ref{charKac} are independent of $\cU$ and
$\cV$.
\end{proof}
\subsection{Extensions of l.c.\ groups} Recall that any normal
$*$-homo\-mor\-phism $\zeta : M_1 \recht \Mh$ of l.c.\ quantum
groups satisfying $\deh \zeta = (\zeta \ot \zeta) \de_1$ generates
two canonical actions: $\mu$ of $(\Mh_1,\hat\de_1)$ on $M$ and
$\theta$ of $(\Mh_1,\hat\de_1^{op})$ on $M$ (see \cite{VV1}). On a
formal level, this means that $\zeta$ gives rise to a dual
morphism $\tilde\zeta: M \recht \Mh_1$ and $\mu$ should be thought
of as $\mu = (\tilde\zeta \ot \id)\de$, while $\te$ should be
thought of as $\te = (\tilde\zeta \ot \id)\deop$.
\begin{definition} \label{32}
Let $G_i\ (i=1,2)$ be l.c.\ groups and let $(M,\de)$ be a l.c.\
quantum group. We call
$$
(L^\infty(G_2),\de_2)\overset{\eta}{\recht} (M,\de)
\overset{\zeta}{\recht} (\mathcal{L}(G_1), \hat\de_1)
$$
a short exact sequence, if
$$
\eta:L^\infty(G_2) \recht M \quad\text{and}\quad \tilde\zeta:
L^\infty(G_1) \recht \Mh
$$
are normal, faithful $*$-homomorphisms satisfying
$$
\de \eta = (\eta \ot \eta) \de_2 \quad\text{and}\quad \deh \zeta =
(\zeta \ot \zeta) \de_1,
$$
and if $\eta(L^\infty(G_2)) = M^\theta$, where $\theta$ is the
canonical action of $(\mathcal{L}(G_1),\hat\de_1^{op})$ on $M$
generated by the morphism $\zeta$. Then we call $(M,\de)$ an
extension of $G_2$ by $\hat{G}_1$.
\end{definition}
The exactness of the sequence in the first, third and second term
is reflected respectively by the faithfulness of $\eta$ and
$\zeta$ and by the formula $\eta(L^\infty(G_2))=M^\theta$.

Given a cocycle matched pair of l.c.\ groups, one can check that
their cocycle bicrossed product is an extension in the sense of
Definition~\ref{32}. Moreover, it belongs to a special class of
extensions, called {\it cleft} extensions (\cite{VV1},
Theorem~2.8). This theorem and \cite{BSV} also show that, whenever
$(M,\de)$ is a cleft extension of $G_2$ by $\hat{G}_1$, then
$(G_1,G_2\subset G)$ is a cocycle matched pair and $(M,\de)$ is
isomorphic to their cocycle bicrossed product.

By definition, two extensions
$$
(L^\infty(G_2),\de_2) \overset{\eta_a}{\recht} (M_a,\de_a)
\overset{\zeta_a}{\recht} (\mathcal{L}(G_1),\hat\de_1)$$ and
$$(L^\infty(G_2),\de_2) \overset{\eta_b}{\recht}
(M_b,\de_b) \overset{\zeta_b}{\recht} (\mathcal{L}(G_1),\hat\de_1)
$$
are called isomorphic, if there is an isomorphism $\pi :
(M_a,\de_a) \recht (M_b,\de_b)$ of l.c.\ quantum groups satisfying
$\pi \eta_a = \eta_b$ and $\hat\pi \zeta_a = \zeta_b$, where
$\hat\pi$ is the canonical isomorphism of $({\Mh}_a,\hat\de_a)$
onto $({\Mh}_b,\hat\de_b)$ associated with $\pi$.

Given a matched pair $(G_1,G_2\subset G)$ of l.c.\ groups, any
couple of cocycles $(\cU,\cV)$ satisfying (\ref{cocU}) generates
as above a cleft extension
$$
(L^\infty(G_2),\de_2 \overset{\eta}{\recht} (M,\de)
\overset{\zeta}{\recht} (\mathcal{L}(G_1),\hat\de_1).
$$
The extensions given by two pairs of cocycles $(\cU_a,\cV_a)$ and
$(\cU_b,\cV_b)$, are isomorphic if and only if there exists a
measurable map $\mathcal{R}$ from $G_1 \times G_2$ to
$\mathbb{T}$, satisfying
$$
\cU_b(g,h,s) = \cU_a(g,h,s) \; \mathcal{R}(h,s) \;
\mathcal{R}(g,\al_h(s)) \; \bar{\mathcal{R}(gh,s)},
$$
$$
\cV_b(g,s,t) = \cV_a(g,s,t) \; \mathcal{R}(g,s) \;
\mathcal{R}(\be_s(g),t) \; \bar{\mathcal{R}}(g,ts)
$$
almost everywhere. Such pairs $(\cU_a,\cV_a)$ and $(\cU_b,\cV_b)$
will be called {\it cohomologous}. The set of equivalence classes
of cohomologous pairs of cocycles $(\cU,\cV)$ satisfying (2.2),
exactly corresponds to the set $\Gamma$ of classes of isomorphic
extensions associated with $(G_1,G_2\subset G)$.

The set $\Gamma$ can be given the structure of an abelian group by
defining
$$
\pi(\cU_a,\cV_a) \; \cdot \; \pi(\cU_b,\cV_b) =
\pi(\cU_a \cU_b,\cV_a \cV_b),
$$
where $\pi(\cU,\cV)$ denotes the equivalence class containing the
pair $(\cU,\cV)$. The group $\Gamma$ is called {\it the group of
extensions} of $(L^\infty(G_2),\de_2)$ by
$(\mathcal{L}(G_1),\hat\de_1)$ associated with the matched pair of
l.c.\ groups $(G_1,G_2\subset G)$. The unit of this group
corresponds to the class of cocycles cohomologous to trivial. The
corresponding extension is called {\it split extension}; all other
extensions are called {\it non-trivial extensions}.
\end{section}
\begin{section}{Regularity properties}
\subsection{General result} As we have seen, one can associate a
multiplicative unitary with any l.c. quantum group. Vice versa,
Baaj and Skandalis constructed in \cite{B}, \cite{BS} a couple of
Hopf $C^*$-algebras in duality out of a given multiplicative
unitary verifying certain regularity conditions. In order to
discuss these conditions for multiplicative unitaries coming from
the bicrossed product construction, let us present the
$C^*$-algebraic version of the split extension (i.e.,\ with
trivial cocycles) and its dual.

Let us associate with a multiplicative unitary $W$ acting on a
Hilbert space $H$, three natural algebras:
$$
S=[(\omega\otimes\id)(W)\vert\omega\in B(H)_{*}],\ \hat
S=[(id\otimes\omega)(W)\vert\omega\in B(H)_{*}],
$$
$$
\text{and}\ [\mathcal{C}(W)]=[(id\otimes\omega)(\Sigma
W)\vert\omega\in B(H)_{*}],
$$
where $[\cdot]$ denotes norm closure. They are not $C^*$-algebras,
in general.
\begin{definition}\label{reg} (see \cite{B}, \cite{BS})
A  multiplicative unitary $W$ is called {\it regular} if
$[\mathcal{C}(W)]=\mathcal{K}(H)$, and {\it semi-regular} if
$[\mathcal{C}(W)]$ contains $\mathcal{K}(H)$, the algebra of all
compact operators on $H$.
\end{definition}
All $W$ associated with Kac algebras, are regular \cite{BS}, and
all $W$ associated with known "isolated" examples of l.c. quantum
groups (see Introduction), are semi-regular; so it was quite a
surprise to find an example of a l.c. quantum group whose
multiplicative unitary is non-semi-regular (see below). If $W$ is
associated with a l.c. quantum group $(M,\de)$, then all
$[\mathcal{C}(W)],\ S\ \text{and}\ \hat S$ are a $C^*$-algebras,
and the comultiplications $\de$ and $\deh$ restrict nicely to
morphisms $S\to M(S\ot S)$ and $\hat S\to M(\hat S\ot\hat S)$
respectively, where $M(A)$ is the multiplier $C^*$-algebra of a
$C^*$-algebra $A$.

If $(M,\de)$ is given by a bicrossed product construction out of a
matched pair $(G_1,G_2\subset G)$ of l.c. groups, then one can
identify $S$ and $\hat S$ with $C_0(G_2 \backslash G)
\rtimes_{\tilde{\be}} G_2$ and $G_1 \kruisje{\tilde{\al}}
C_0(G/G_1)$ respectively, where $\tilde\al_g$ and $\tilde\be_s$
are the canonical continuous actions of $G_1$ on $G/G_1$ of $G_2$
on $G_2 \backslash G$ respectively ($C_0(X)$ denotes the
$C^*$-algebra of continuous functions vanishing at infinity on a
l.c. topological space $X$) \cite{BSV}. One can check that the
measurable mutual actions $\al_g$ and $\be_s$ of $G_1$ and $G_2$
are the restrictions of the above canonical continuous actions
$\tilde\al_g$ and $\tilde\be_s$ (topologies on $G_1$ and $G_2
\backslash G$ and, respectively, on $G_2$ and $G/G_1$, are in
general different).

Now we can formulate the main result of \cite{BSV}:
\begin{theorem}\label{sreg}
The multiplicative unitary $W$ of the bicrossed product l.c.
quantum group $(M,\de)$ is regular if and only if the map
$$
\te : G_1\times G_2\to G : \te(g,s)=gs
$$
is a homeomorphism of $G_1\times G_2$ onto $G$. $W$ is
semi-regular if and only if $\te$ is a homeomorphism of $G_1\times
G_2$ onto an open subset of $G$ of full measure.
\end{theorem}

We will call the corresponding matched pairs of l.c. groups
regular, semi-regular and non-semi-regular, respectively. To get
an example of a regular matched pair, it suffices to take $G$
discrete or in the form of a semi-direct product of $G_1$ and
$G_2$, with closed subgroups $G_1$ and $G_2$. In both cases the
bicrossed product l.c. quantum group is a Kac algebra, due to
Corollaries \ref{411} and \ref{410}, respectively.

\begin{remark} Woronowicz \cite{W5} constructed a couple of Hopf
$C^*$-algebras in duality out of a given multiplicative unitary
under certain alternative conditions of {\it manageability}. All
multiplicative unitaries associated with l.c. quantum groups, are
manageable.
\end{remark}
\subsection{Matched pairs of Lie groups \cite{VV2}} We consider Lie
groups and Lie algebras over the field $\mathbb{R}$ or
$\mathbb{C}$. A matched pair of {\it Lie groups} (i.e.,\ when $G$
is a Lie group), is always semi-regular:

\begin{proposition} \label{LieOK}
If, in Definition~\ref{mpg}, $G$ is a Lie group, then the map
$\te$ has an open range $\Omega$ and is a diffeomorphism of $G_1
\times G_2$ onto $\Omega$, where $G_1$ and $G_2$ are Lie groups
under the identification with closed subgroups of $G$.
\end{proposition}

The infinitesimal form of the last notion is as follows (see
\cite{Majbook}):

\begin{definition} \label{def2.2}
We call $(\lieg_1,\lieg_2)$ a matched pair of Lie algebras, if
there exists a Lie algebra $\lieg$ with Lie subalgebras $\lieg_1$
and $\lieg_2$ such that $\lieg = \lieg_1 \oplus \lieg_2$ as vector
spaces.
\end{definition}

These conditions are equivalent \cite{Majbook} to the existence of
a left action $\triangleright: \lieg_2\otimes \lieg_1\to \lieg_1$
and a right action $\triangleleft: \lieg_2\otimes \lieg_1\to
\lieg_2$, so that $\lieg_1$ is a left $\lieg_2$-module and
$\lieg_2$ is a right $\lieg_1$-module and
\begin{enumerate}
\item $x\triangleright[a,b]=[x\triangleright a,b]+[a,x\triangleright b]+
(x\triangleleft a)\triangleright b-(x\triangleleft
b)\triangleright a, $
\item $[x,y]\triangleleft a=[x,y\triangleleft a]+[x\triangleleft a,y]+
x\triangleleft (y\triangleright a)-y\triangleleft (x\triangleright
a), $
\end{enumerate}
for all $a,b\in \lieg_1,\ x,y\in \lieg_2$. Then, for $\lieg =
\lieg_1 \oplus \lieg_2$ we have:
$$
[a\oplus x,b\oplus y]=([a,b]+x\triangleright b - y\triangleright
a)\oplus ([x,y]+x\triangleleft b-y\triangleleft a).
$$

Two matched pairs, $(\lieg_1,\lieg_2)$ and $(\lieg'_1,\lieg'_2)$,
are called isomorphic if there is an isomorphism of the
corresponding Lie algebras $\lieg$ and $\lieg'$ sending $\lieg_i$
onto $\lieg'_i\ (i=1,2)$.

Let us explain the relation between the two notions of a matched
pair.

\begin{proposition} \label{prop2.3}
Let $(G_1,G_2\subset G)$ be a matched pair of Lie groups. If
$\lieg$ denotes the Lie algebra of $G$, and if $\lieg_1$, resp.\
$\lieg_2$, are the Lie subalgebras corresponding to the closed
subgroups $i(G_1)$, resp.\ $j(G_2)$, then $(\lieg_1,\lieg_2)$ is a
matched pair of Lie algebras.
\end{proposition}
\begin{proof}
The fact that $\lieg = \lieg_1 \oplus \lieg_2$ as vector spaces
follows from the fact that $\te$ is a diffeomorphism in the
neighborhood of the unit element.
\end{proof}

The converse problem, to construct a matched pair of Lie groups
from a given matched pair $(\lieg_1,\lieg_2)$ of Lie algebras, is
much more delicate. Let us start with

\begin{proposition} \label{propnieuw}
Every matched pair of Lie algebras $\lieg_1=\lieg_2=k$ can be
exponentiated to a matched pair of Lie groups $(G_1,G_2\subset
G)$, where $G_1,G_2$ are either $(k,+)$ or $(k\setminus
\{0\},\cdot)$, where $k=\mathbb{C}\ \text{or}\ \mathbb{R}$.
\end{proposition}
\begin{proof}
a) The only two-dimensional complex Lie algebras are the abelian
one and the one with generators $X,Y$ and relation $[X,Y] =Y$. If
$\lieg$ is abelian, the mutual actions of $\lieg_1$ and $\lieg_2$
on each other are trivial and exponentiation is a direct sum.

If $\lieg$ is generated by $[X,Y]=Y$, we either have that
$\lieg_1$ or $\lieg_2$ equals to $\mathbb{C} Y$, then one of the
actions is trivial and $G$ can be constructed as semi-direct
product of the connected, simply connected Lie groups of $\lieg_1$
and $\lieg_2$, or both $\lieg_1,\ \lieg_2\neq\mathbb{C} Y$. In the
latter case, there is, up to isomorphism, only one possibility,
namely $\lieg_1 = \mathbb{C} X$, $\lieg_2 = \mathbb{C} (X+Y)$.
Define on $\mathbb{C} \setminus \{0\} \times \mathbb{C}$ the Lie
group with product
$$
(t,s)(t',s')=(tt',s+ts') \; .
$$
Define $G_1=G_2=\mathbb{C} \setminus \{0\}$ with embeddings $i(g)
= (g,0)$ and $j(s) = (s,s-1)$, we indeed get a matched pair of
complex Lie groups with mutual actions
\begin{equation} \label{exp11}
\al_g(s) = g(s-1) + 1 \; , \quad \be_s(g) = \frac{sg}{g(s-1) + 1}
\; .
\end{equation}
b) The real case is completely analogous.
\end{proof}

\begin{remark} \label{remnieuw} The natural idea to take the
connected, simply connected Lie group $G$ of the corresponding
$\lieg$ and its unique connected, closed subgroups $G_1$ and $G_2$
with tangent Lie algebras $\lieg_1$ and $\lieg_2$, respectively,
fails because such $(G_1,G_2)$ is not necessarily a matched pair
even if $\dim \lieg_1=\dim \lieg_2=1$. Indeed, if $k=\mathbb{C}$,
the connected simply connected Lie group $G$ of $\lieg$ consists
of all pairs $(t,s)$ with $t,s\in \mathbb{C}$ and the product
$$
(t,s)(t',s')=(t+t',s+\exp(t)s'),
$$
and its closed subgroups $G_1$ and $G_2$ corresponding to the
decomposition $\lieg = \mathbb{C} X \oplus \mathbb{C} (X+Y)$ above
consist respectively of all pairs of the form $(g,0)$ and
$(s,\exp(s)-1)$ with $g,s\in \mathbb{C}$. These groups do not form
a matched pair because $G_1\cap G_2=\{(2\pi in,0)\vert n\in
\mathbb{Z}\}=Z(G)$. So, it is crucial not to take $G$ simply
connected.

Taking $g,t,s$ above real, we come to the example of a matched
pair of real Lie groups from \cite{VV1}, Section~5.3 (see also
\cite{Skand} and \cite{SV2}). Here $\lieg$ is a real Lie algebra
generated by $X$ and $Y$ subject to the relation $[X,Y]=Y$ and one
considers the decomposition $\lieg=\mathbb{R}
X\oplus\mathbb{R}(X+Y)$. Then, to get a matched pair of Lie
groups, we consider $G$ as the variety $\mathbb{R} \setminus \{0\}
\times \mathbb{R}$ with the product
$$
(s,x)(t,y) = (st,x+sy)
$$
and embed $G_1=G_2=\mathbb{R} \setminus \{0\}$ by the formulas
$i(g) = (g,0)$ and $j(s) = (s,s-1)$. Remark that here, it is
impossible to take the connected component of the unity of the
group of affine transformations of the real line as $G$, because
it is easy to see that for its closed subgroups $G_1$ and $G_2$
corresponding to the above mentioned subalgebras, the set $G_1G_2$
is not dense in $G$. The corresponding multiplicative unitary is
semi-regular, but not regular.
\end{remark}

It is even possible that for a given matched pair of Lie algebras,
$G_1 \cap G_2 \neq \{e\}$ for {\em any} corresponding pair of Lie
groups, i.e., such a matched pair of Lie algebras cannot be
exponentiated to a matched pair of Lie groups.

\begin{example} \label{ex2.4}
Consider a family of complex Lie algebras
$\lieg=\text{span}\{X,Y,Z\}$ with $[X,Y]=Y,\ [X,Z]=\al Z,\
[Y,Z]=0$, where $\al\in \mathbb{C}\setminus \{0\}$, and the
decomposition $\lieg=\text{span}\{X,Y\}\oplus\mathbb{C} (X+\al
Z)$. The corresponding connected simply connected complex Lie
group $H$ consists of all triples $(t,u,v)$ with $t,u,v\in
\mathbb{C}$ and the product
$$(t,u,v)(t',u',v')=(t+t',u+\exp(t)u',v+\exp(\al t)v'),$$
and its closed subgroups $H_1$ and $H_2$ corresponding to the
decomposition above consist respectively of all triples of the
form $(t,u,0)$ and $(s,0,\exp(\al s)-1)$ with $t,u,s\in \C$. These
groups do not form a matched pair because $H_1\cap
H_2=\{({\frac{2\pi in}{\al}},0,0)\vert n\in \mathbb{Z}\}$.

We claim that, if $1/\al \not\in \mathbb{Z}$ and if $G$ is any
complex Lie group with Lie algebra $\lieg$, such that $G_1,G_2$
are closed subgroups of $G$ with tangent Lie algebras $\lieg_1$,
resp.\ $\lieg_2$, then $G_1 \cap G_2 \neq \{e\}$. Indeed, since
the Lie group $H$ is connected and simply connected, the connected
component $G^{(e)}$ of $e$ in $G$ can be identified with the
quotient of $H$ by a discrete central subgroup. If $\al \not\in
\mathbb{Q}$, the center of $H$ is trivial, so that we can identify
$G^{(e)}$ and $H$. Under this identification, the connected
components of $e$ in $G_1,G_2$ agree with $H_1,H_2$. Because $H_1
\cap H_2 \neq \{e\}$, our claim follows. If $\al =\frac{m}{n}$ for
$m,n \in \mathbb{Z} \setminus \{0\}$ mutually prime, the center of
$H$ consists of the elements $\{(2\pi n N,0,0) \mid N \in
\mathbb{Z} \}$. Hence, the different possible quotients of $H$ are
labeled by $N \in \mathbb{Z}$ and are given by the triples
$(a,u,v) \in \mathbb{C}^3$, $a \neq 0$ and the product
\begin{equation} \label{group-law}
(a,u,v) (a',u',v') = (aa',u+a^{nN}u',v+a^{mN}v') \; .
\end{equation}
The closed subgroups corresponding to $\lieg_1$ and $\lieg_2$ are
given by $(a,u,0)$ and $(b,0,b^{mN} - 1)$ with $a,b,u \in
\mathbb{C}$ and $a,b \neq 0$. The intersection of both subgroups
is non-trivial whenever $mN \neq \pm 1$. This proves our claim.

Considering now the complex Lie algebras above as real Lie
algebras with generators $X,iX,Y,iY,Z,iZ$ and the decomposition
above as a decomposition of real Lie algebras, we get a matched
pair of real Lie algebras which cannot be exponentiated to a
matched pair of real Lie groups.

In the remaining case $\al=1/n$ with $n \in \mathbb{Z} \setminus
\{0\}$, we can consider the Lie group $G$ defined by
Equation~\eqref{group-law} with $m=N=1$. Consider $G_1 =
\mathbb{C} \setminus \{0\} \times \mathbb{C}$ with $(a,u)(a',u') =
(aa',u+a^n u')$ and $G_2 = \mathbb{C} \setminus \{0\}$. Writing
$i(a,u) = (a,u,0)$ and $j(v) = (v,0,v-1)$, we get a matched pair
of Lie groups with mutual actions
$$\al_{(a,u)}(v) = a(v-1)+1 \quad\text{and}\quad \be_v(a,u) = \bigl(
\frac{va}{a(v-1)+1}, \frac{u}{(a(v-1)+1)^n} \bigr)\; .$$
\end{example}

However, any matched pair of real Lie algebras when one of them
has dimension 1 and the other at most 2, can be exponentiated.

\subsection{Cocycle matched pairs of Lie groups \cite{VV2}}
The usage of 2-cocycles gives much more concrete examples of l.c.\
quantum groups. We first explain the infinitesimal picture, i.e.
how 2-cocycles for matched pairs of Lie algebras look like, how
they are related to extensions and then discuss the problem of
exponentiation.

Recall that a Lie bialgebra, according to V.G. Drinfeld, is a Lie
algebra $\lieg$ equipped with a Lie bracket $[\cdot,\cdot ]$ and a
Lie cobracket $\delta$, i.e., a linear map
$\delta:\lieg\to\lieg\otimes \lieg$ satisfying the
co-anticommutativity and the co-Jacobi identity, that is:
$$
(\id - \tau)\delta = 0 \; , \quad  (\id + \zeta +
\zeta^2)(\id\otimes\delta)\delta=0,
$$
where
$$
\tau(u\otimes v)=v\otimes u, \quad \zeta(u\otimes v\otimes
w)=v\otimes w\otimes u \quad \ (\text{for all}\quad u,v,w\in
\lieg)
$$
are the flip maps, and $[\cdot,\cdot],\ \delta$ are compatible in
the following sense:
$$
\delta[u,v]=[u,v_{[1]}]\otimes v_{[2]}+v_{[1]}\otimes
[u,v_{[2]}]+[u_{[1]},v]\otimes u_{[2]} +u_{[1]}\otimes
[u_{[2]},v].
$$
Any Lie algebra (respectively, Lie coalgebra, i.e., vector space
dual to a Lie algebra) is a Lie bialgebra with zero Lie cobracket
(respectively, zero Lie bracket). The definition of a morphism of
Lie bialgebras is obvious.

Given a pair of Lie algebras $(\lieg_1,\lieg_2)$, let look for a
Lie bialgebra $\lieg$ such that
$$
\lieg_2^* \longrightarrow \lieg \longrightarrow \lieg_1
$$
is a short exact sequence in the category of Lie bialgebras. This
means precisely that $\lieg$ has a sub-bialgebra with trivial
bracket, which is an ideal and such that the quotient is a Lie
bialgebra with trivial cobracket.

The theory of extensions in this framework has been developed in
\cite{Mas2} and is similar to the theory of extensions of l.c.\
groups that we discussed above. Namely, for the existence of an
extension $\lieg$, it is necessary and sufficient that $(\lieg_1,
\lieg_2)$ form a matched pair, and all extensions are bicrossed
products with cocycles. We consider this theory as an
infinitesimal version of the theory of extensions of Lie groups.

As we remember, for any matched pair of Lie algebras $(\lieg_1,
\lieg_2)$, there are mutual actions $\triangleright:
\lieg_2\otimes \lieg_1\to \lieg_1$ and $\triangleleft:
\lieg_2\otimes \lieg_1\to \lieg_2$, compatible in a way explained
above and such that, for all $a,b\in \lieg_1,\ x,y\in \lieg_2$, we
have
$$
[a\oplus x,b\oplus y]=([a,b]+x\triangleright b - y\triangleright
a)\oplus ([x,y]+x\triangleleft b-y\triangleleft a).
$$

The definition of a pair of 2-cocycles on such a matched pair is
given in \cite{Majbook}, \cite{Mas2}. For our needs, it suffices
to understand that these 2-cocycles are linear maps
$$
\cU : \lieg_1 \wedge \lieg_1 \recht \lieg_2^* \; , \qquad \cV :
\lieg_2 \wedge \lieg_2 \recht \lieg_1^* \;
$$
verifying certain 2-cocycle equations and compatibility equations
that are infinitesimal forms of equations (\ref{cocU}). The link
between 2-cocycles on matched pairs of real Lie groups and those
of real Lie algebras is given by
\begin{proposition} \label{cocycles}
Let $(G_1,G_2)$ be a matched pair of real Lie groups equipped with
cocycles $\cU$ and $\cV$, which are differentiable around the unit
elements, and let $(\lieg_1, \lieg_2)$ be the corresponding
matched pair of real Lie algebras. Defining
$$
\langle \cU(X,Y), A \rangle = -i (X_e \ot Y_e \ot A_e - Y_e \ot
X_e \ot A_e)(\cU) \quad\text{and} $$ $$ \langle \cV(A,B), X
\rangle = -i (A_e \ot B_e \ot X_e - B_e \ot A_e \ot X_e)(\cV) \; ,
$$
for $X,Y \in \lieg_1$ and $A,B \in \lieg_2$, we get a pair of
cocycles on $(\lieg_1, \lieg_2)$.
\end{proposition}
Here $\langle \cdot , \cdot \rangle$ is the duality between
$\lieg_i$ and $\lieg_i^*$, and $X_e,Y_e,A_e,B_e$ denote the
partial derivatives at $e$ in the direction of the corresponding
generator. The factor $-i$ appears because for Lie groups $\cU$
and $\cV$ take values in $\mathbb{T}$, and for real Lie algebras
we consider 2-cocycles as real linear maps.

It is clear that these 2-cocycles $\cU$ form a real vector space.
If one of the Lie groups is 1-dimensional, the corresponding
2-cocycle is trivial, since it can be regarded as an
antisymmetric, bilinear form on the corresponding Lie algebra.

We defined above the group of extensions for a matched pair of
l.c. groups $(G_1,G_2)$ using the notion of cohomologous
2-cocycles. The same can be done for a matched pair of Lie
algebras \cite{Mas2}. Two cocycles $\cU_1$ and $\cU_2$ are called
\emph{cohomologous} if $\cU_1-\cU_2$ is cohomologous to trivial.
The quotient space of 2-cocycles modulo 2-cocycles cohomologous to
trivial, with addition as the group operation, is called the
\emph{group of extensions} of the matched pair $(\lieg_1,
\lieg_2)$.

In particular, the group of extensions for a matched pair of real
Lie algebras of dimensions $1$ and $2$, is either trivial or
$\mathbb{R}$. To exponentiate these 2-cocycles, i.e., to construct
the corresponding 2-cocycles on the level of Lie groups, is much
more difficult. However, a complete classification of cocycle
matched pairs of real Lie groups with $dim G\leq 3$ and at most
$2$ was obtained in \cite{VV2}; the corresponding groups of
extensions are either trivial, or $\mathbb{R}$, or $\mathbb{Z}$.

A series of examples of l.c. quantum groups was constructed in
\cite{CKP1} out of cocycle matched pairs of real 2-dimensional Lie
algebras by their exponentiation.
\subsection{Non-semi-regular l.c. quantum groups \cite{BSV}} Let
$\cA\neq 0$ be a l.c. ring with additive Haar measure $\nu$, and
let $\cA^\times$ be its group of invertible elements. Then
$\cA^\times$ is a l.c. group by considering it as a closed
subspace $\{(a,b)\vert ab=ba=1\}$ of $\cA\times\cA$, let
$\nu^\times$ be its Haar measure. Consider the group
$G=\{(a,x)\vert a\in\cA^\times, x\in\cA\}$ with multiplication
$(a,x)(b,y)=(ab,x+ay)$, and its closed subgroups
$G_1=\{(a,a-1)\vert a\in\cA^\times\}$, $G_2=\{(b,0)\vert
b\in\cA^\times\}$. The Haar measure of $G$ is the product of $\nu$
and $\nu^\times$. Then $G_1G_2=\{(a,x)\in \cA^\times\times\cA\vert
x+1\in\cA^\times\}$. Hence, $(G_1,G_2)$ is a matched pair if and
only if $A^\times$ has complement of measure $\nu$ zero in $\cA$.
Then one can prove
\begin{proposition}\label{A}
The above matched pair is not regular. It is semi-regular if and
only if $\cA^\times$ is open in $\cA$. If $(M,\de)$ is the
corresponding bicrossed product l.c. quantum group, then the dual
$(\Mh,\deh)$ is isomorphic to the opposite quantum group
$(M,\de^{op})$ and $M\widetilde{=}\hat M\widetilde{=}
B(L^2(\cA^\times,\nu^\times))$.
\end{proposition}

\begin{example}\label{nsr}
In order to construct an example of a non-semi-regular matched
pair, let us choose the second countable ring $\cA$ in the
following way. Let $\cP$ be a set of prime numbers such that
$$
\Sigma_{p\in{\cP}}\frac{1}{p}<\infty.
$$
Define the restricted Cartesian product of l.c. fields
$\mathbb{Q}_{p}$ of $p$-adic rational numbers relatively to the
compact open subrings $\mathbb{Z}_{p}$ of the corresponding
$p$-adic integers:
$$
\cA=\Pi'_{p\in\mathcal{P}}(\mathbb{Q}_{p},\mathbb{Z}_{p}).
$$
It consists of sequences $(x_p)\in\mathbb{Q}_{p}$ which eventually
belong to $\mathbb{Z}_{p}$. Equipped with the usual l.c. topology,
$\cA$ becomes a l.c. ring such that $\cA^\times$ has complement of
measure zero and empty interior in $\cA$. The fact that
$\cA^\times$ has complement of measure zero in $\cA$ follows from
the Borel-Cantelli lemma: normalizing the Haar measure on
$\mathbb{Q}_{p}$ in such a way that $\mathbb{Z}_{p}$ has measure
one, observe that $\mathbb{Z}_{p}\backslash\mathbb{Z}_{p}^\times$
has measure $\frac{1}{p}$, which assumed to be summable over
$p\in\cP$.

Proposition \ref{A} shows that the corresponding multiplicative
unitary in not semi-regular. It is known that the $C^*$-algebra
$S=\cA^\times\ltimes C_0(\cA)$ is not of type I \cite{BC}.
\end{example}
\end{section}
\begin{section}{L.c. quantum groups whose von Neumann algebras are
factors}
\subsection{Motivation} Let $\mathcal{L}(G)$ be the von Neumann algebra
generated by the operators $\lambda_g,\ g\in G$ of the left
regular representation of a l.c. group $G$ or, equivalently, by
the operators of the form $L_f=\int_G f(g)\lambda_g d\mu_G$, where
$f\in L^1(G,\mu_G)$. All these operators act on the Hilbert space
$L^2(G,\mu_G)$, where $\mu_G$ is a left Haar measure on $G$. It is
equipped with the canonical n.s.f. weight defined by
$\varphi(L_f)=f(e)$, for all functions from $L_1(G,\mu_G)$
continuous in the neutral element $e\in G$. This weight is a trace
if and only if $G$ is unimodular.

If $G\neq \{e\}$ is compact, then $\mathcal{L}(G)$ is a direct sum
of finite dimensional full matrix algebras and cannot be a factor
(i.e., a von Neumann algebra with trivial center). Indeed, in this
case one of the summands corresponding to the trivial
representation of $G$ must be of dimension 1. But there exist
non-compact groups $G$ with $\mathcal{L}(G)$ a factor of any type,
in the sense of Murray-von Neumann classification, except for a
full matrix algebra of finite dimension: I$_\infty$, II$_1$,
II$_\infty$, III$_\lambda$, where $\lambda\in [0,1]$ (see, for
example, \cite{JS}; the standard reference on the classification
of type III factors and their invariants is \cite{C}).

a) For the group $G=\mathbb{R}\rtimes\mathbb{R}^\times$ -
semi-direct product with natural action of $\mathbb{R}^\times$ on
$\mathbb{R}$ by multiplication, one has $\mathcal{L}(G)\simeq
B(L^2(\mathbb{R}^\times)$ - the type I$_\infty$ factor.

b) Let $G$ be a discrete countable group with all nontrivial
conjugacy classes infinite, for instance, $\mathbb{F}_n$ - the
free group with $n\geq 2$ generators or $S_\infty$ - the group of
permutations of the set of natural numbers such that any
individual permutation permutes only finitely many numbers. Then
$\mathcal{L}(G)$ equipped with the canonical finite trace
$\varphi$ is a type II$_1$ factor (see, for example, \cite{JS}).

c) The Cartesian product of groups in a) and b) gives an example
of $G$ such that $\mathcal{L}(G)$ is a type II$_\infty$ factor.

d) R. Godement showed that for $G=\mathbb{R}^2\rtimes
GL_2(\mathbb{Q})$  - semi-direct product with natural action of
$GL_2(\mathbb{Q})$ on $\mathbb{R}^2$, $\mathcal{L}(G)$ is a
non-hyperfinite type III$_1$ factor. Then C. Sutherland derived
from this a series of examples of groups for which
$\mathcal{L}(G)$ is a non-hyperfinite type III$_\lambda$ factor,
for all $\lambda\in [0,1)$, and A. Connes constructed similar
examples with hyperfinite factors, for $\lambda\in (0,1)$ (for all
these results see \cite{Su}). Later on, for all $\lambda\in
[0,1]$, examples of groups $G$ coming from number theory with
$\mathcal{L}(G)$ a hyperfinite type III$_\lambda$ factor, were
constructed in \cite{BZ}. Moreover, one can construct such
examples with specific properties of the invariant $T$ (see
\cite{Su}, \cite{BZ}).

These results and certain freedom given by the bicrossed product
construction in producing examples of l.c. quantum groups, allows
us to ask if it is possible to construct such examples of
$(M,\de)$ that both von Neumann algebras, $M$ and $\hat M$, are
factors of all prescribed types, except for a full matrix algebra
of finite dimension. Such quantum groups would be "as far as
possible" from usual groups since, by definition, a factor is a
von Neumann algebra with trivial center. We present below some
partial results in this direction.

Since the cocycle bicrossed product construction gives $M=G_1
\kruisje{\al,\cU} L^{\infty}(G_2)$ and $\hat M=G_2
\kruisje{\al,\cU} L^{\infty}(G_1)$, let us look when the centers
of these algebras equal to $\mathbb C 1$. If $\cU$ and $\cV$ are
trivial, one can use the result from \cite{Sauv}. Recall that a
borelian action $\te$ of a l.c. group $\Gamma$ on a l.c.
topological space $X$ equipped with a borelian measure $\mu$, is
called {\it free} if the stabilizer $S_x:=\{\gamma\in \Gamma\vert
gx=x\}$ equals to $\{e\}$ for $\mu$-almost all $x\in X$. For such
an action, \cite{Sauv} Corollary 2.3 says that the center of
$\Gamma\kruisje{\te} L^{\infty}(X,\mu)$ is a factor if and only if
$\te$ is also {\it ergodic} (i.e., the set of $\te$-fixed points
is of measure zero). Thus, it would be interesting to find matched
pairs $(G_1,G_2)$ with both actions $\al$ and $\be$ being free and
ergodic (remark immediately that if $G_1$ (resp., $G_2$) is
discrete, the relations $\al_g(e)=e,\ \be_s(e)=e,\ \forall g\in
G_1,s\in G_2$ show that $\be$ (resp., $\al$) is not free).

One instance of this kind is given by matched pairs with $G_2= g_0
G_1g_0^{-1}$, where $G_0$ is an element of $G$, i.e., $G_1$ and
$G_2$ are conjugated. Then \cite{BSV} Prop. 3.13 says that both
actions $\al$ and $\be$ are isomorphic to the action of $G_2$ on
itself by translations which is clearly free and ergodic.
Moreover, in this case both $M$ and $\hat M$ are isomorphic to
$B(L^2(G_2))$, i.e., are type I$_\infty$ factors.

In particular, for the matched pares introduced before Proposition
\ref{A}, $G_1$ and $G_2$ are conjugated by the element $(-1,-1)$.
Another example (\cite{BSV}, Ex. 4.6) is as follows. Let again
$\cA\neq 0$ be l.c. ring such that $A^\times$ has complement of
measure $\nu$ zero in $\cA$. Take $G=GL_2(\cA)$ and
$G_1=\{(a_{ij})\}$ with $a_{11}\in\cA^\times,\ a_{12}\in\cA,\
a_{21}=0,\ a_{22}=1$,\ $G_2=\{(a_{ij})\}$ with $a_{11}=1,\
a_{12}=0,\ a_{21}\in\cA,\ a_{22}\in\cA^\times$. Clearly, $G_1$ and
$G_2$ are conjugated.

\subsection{L.c. quantum groups whose algebras are type II and
type III factors \cite{F}}
Using the techniques close to that of \cite{W1}, one can show that
a l.c. quantum group $(M,\de)$ such that $M$ is a type II$_1$
factor, is necessarily a compact quantum group in the sense of
\cite{W1}, so that $\hat M$ is necessarily a direct sum of finite
dimensional full matrix algebras and cannot be a factor (one of
the summands is generated by the counit and must be of dimension
1). Thus, the case when $M$ or $\hat M$ is a type II$_1$ factor,
must be excluded from consideration.

The following examples of l.c. quantum groups whose algebras are
factors, are inspired by Example \ref{nsr}. We use again second
countable l.c. ring of the form
$$
\cA=\Pi'_{n\in\mathbb{N}}(\mathbb{Q}_{p_n},\mathbb{Z}_{p_n}),
$$
where $(p_n)$ is now an infinite sequence of prime numbers without
repetitions, its group of invertible elements
$$
\cA^\times=\Pi'_{n\in\mathbb{N}}(\mathbb{Q}^\times_{p_n},\mathbb{Z}^\times_{p_n}),
$$
the group $G=\{(a,x)\vert a\in\cA^\times, x\in\cA\}$ with
multiplication $(a,x)(b,y)=(ab,x+ay)$, and the closed subgroup
$G_1=\{(a,0)\vert a\in\cA^\times\}$ of $G$. But now we choose the
second closed subgroup $G_2$ of $G$ in another way:
$$
G_2=\{((a_n),(b_n))\in G\vert a_n + b_n p_n=1, \forall
n\in\mathbb{N}\}.
$$
Then, checking the conditions of Theorem \ref{sreg}, one can prove
that $(G_1,G_2\subset G)$ is a matched pair of l.c. groups, which
is semi-regular, but not regular. Thus, one can apply the
bicrossed product construction and to get this way two l.c.
compact quantum groups in duality, whose von Neumann algebras will
be $M=G_1\ltimes L^\infty(G_2)$ and $\hat M=G_2\ltimes
L^\infty(G_1)$. Moreover, one can show that both actions $\al$ and
$\be$ are free and ergodic, so that both $M$ and $\hat M$ are
factors.

In order to determine their types, one constructs their explicit
isomorphisms to so called {\it infinite tensor products of type I
factors}, briefly, ITPFI factors (see, for example, \cite{AW}), of
the form
$$
L=\otimes_{n\in\mathbb{N}} (B(l^2(\mathbb{N})),\varphi_{n}),
$$
where $\varphi_{n}$ is a normal state on the type I factor
$B(l^2(\mathbb{N}))$ of the form $\varphi_{n}(\cdot)=tr(\rho_{n}\
\cdot)$. Here $\rho_{n}$ is a positive trace class operator
characterized by the decreasing sequence of its eigenvalues. In
the case of our $M$, this sequence of eigenvalues is
$\lambda^M_{n,k}=(1-p_n^{-1})p_n^{-k},\ n,k\in \mathbb{N}$, and in
the case of $\hat M$ the corresponding sequence of eigenvalues
$\lambda^{\hat M}_{n,k},\ n,k\in\mathbb{N}$ also can be computed
explicitly. Now we can use the following general result \cite{AW}:

\begin{proposition}\label{aw} Let
$$
L=\otimes_{n\in\mathbb{N}} (M_n,\varphi_n),
$$
be an ITPFI factor, where $M_n$ is a type I$_{n_\nu}$ factor with
$n_\nu\in \mathbb{N}\cup\{\infty\}$, and let $\varphi_n$ be a
normal state on $M_n$ of the form $\varphi_n(\cdot)=tr(\rho_{n}\
\cdot)$. Here $\rho_n$ is a positive trace class operator
characterized by the decreasing sequence of its eigenvalues
$\lambda_{n_i},\ i=1,2,...,n_\nu$. Then:

1. $L$ is of type I if and only if
$$
\Sigma_n\vert 1-\lambda_{n_1}\vert<+\infty.
$$
2. $L$ is of type II if and only if $n_\nu<+\infty$ for all $n$
and
$$
\Sigma_{n,i}\vert
(n_\nu)^{-\frac{1}{2}}-(\lambda_{n_i})^{\frac{1}{2}}\vert^2<+\infty.
$$
3. If exists $\delta>0$ such that for all $n$ we have
$\lambda_{n_1}\geq\delta$, then $L$ is of type III if and only if
$$
\Sigma_{n,i}\inf\{\vert
\frac{\lambda_{n_1}}{\lambda_{n_i}}-1\vert^2,C\}<+\infty,
$$
for some constant $C>0$.
\end{proposition}
Comparing the above mentioned sequences of eigenvalues
$\lambda^M_{n,k}$ and $\lambda^{\hat M}_{n,k},\ n,k\in\mathbb{N}$,
with conditions of Proposition \ref{aw} and using the results of
\cite{BZ}, one gets the following
\begin{theorem} \label{f} \cite{F}
For each infinite sequence $(p_n)$ of prime numbers without
repetitions, the von Neumann algebras $M$ and $\hat M$ are
hyperfinite factors. Then:

(i) Both $M$ and $\hat M$ have the same invariant $T$ (in the
sense of \cite{C}).

(ii) The condition $\Sigma_n\frac{1}{p_n}<+\infty$ is equivalent
to the fact that $\cA^\times$ has complement of measure zero in
$\cA$. In this case, $(M,\Delta)$ is of type
$(I_\infty,II_\infty)$.

(iii) The condition $\Sigma_n\frac{1}{p_n}=+\infty$ is equivalent
to the fact that $\nu(\cA^\times)$ is zero. In this case,
$(M,\Delta)$ is of type $(III,III)$.

Moreover, we have:

a) For each $\lambda\in[0,1]$, there exists such $({p_n})$ that
both $M$ and $\hat M$ are of type III$_\lambda$.

b) For each countable subgroup $K$ of $\mathbb {R}$ and a
countable subset $\Sigma$ of $\mathbb {R}\backslash K$ there
exists such $({p_n})$ that that the corresponding invariant $T$
contains $K$ and does not intersect with $\Sigma$.
\end{theorem}
Remark, that tensor products of $(M,\Delta)$ corresponding to the
case (i) and their duals give examples of l.c. quantum groups with
both factors of type $II_\infty$.
\end{section}
\begin{section}{Some other results}

\subsection{Amenability  \cite{DQV}}

Let $(M,\de)$ be a l.c. quantum group. A state $m\in M^*$ is said
to be a {\it left invariant mean} on $(M,\de)$ if
$$
m((\omega\otimes\id)\de(x))=m(x)\omega(1),
$$
for all $\omega\in M_{*}$ and $x\in M$. Similarly, one defines a
right invariant mean. An invariant mean is both right and left
invariant mean. One can show that a left invariant mean exists on
$(M,\de)$ if and only if an invariant mean exists on it.
\begin{definition} \label{wa}
We call $(M,\de)$ {\it amenable} if there exists a left invariant
mean on it. We say that $(M,\de)$ is {\it coamenable} if
$(\Mh,\deh)$ is {\it amenable}.
\end{definition}
\begin{definition} \label{sa}
We call $(M,\de)$ {\it strongly amenable} if there exists a
bounded counit on $(\hat S,\deh)$, the reduced dual
$C^*$-algebraic l.c. quantum group \cite{KV1}. We say that
$(M,\de)$ is {\it strongly coamenable} if $(\Mh,\deh)$ is {\it
strongly amenable}.
\end{definition}
Strong amenability implies amenability; the two properties
coincide for usual l.c. groups and for discrete l.c. quantum
groups. It is not known if they coincide for general l.c. quantum
groups, even for Kac algebras. The amenability for Kac algebras
was studied in \cite{ES1}.

\begin{remark}\label{Ruan} It was claimed in \cite{ES1} that they
coincide for Kac algebras, but Zh.-J. Ruan found a gap in the
proof. He showed \cite{Ruan} that these properties coincide for
discrete Kac algebras. Later on, E. Blanchard and S. Vaes proved
the same for discrete quantum groups (unpublished), their proof
uses the Powers-St{\o}rmer inequality. One can find another proof
of the last result in \cite{Tom}.
\end{remark}

According to \cite{DQV}, an extension of the form
$$(M_2,\de_2) \recht(M,\de) \recht ({\hat M}_1,
\hat\de_1),$$ where $(M_1,\de_1),\ (M_2,\de_2)$ are two l.c.
quantum groups (see \cite{VV1}) is amenable if and only if both
$(M_1,\de_1)\ \text{and}\ (M_2,\de_2)$ are amenable. The same
holds in strongly amenable case, but only for split extensions.
This allows to construct various examples. Given a matched pair of
amenable l.c. groups, all the corresponding extensions are
amenable, and the split extension is strongly amenable. All
examples in \cite{VV1} are of this kind.

Two concrete examples of non-amenable l.c. quantum groups were
constructed in \cite{DQV} using the techniques of \cite{VV2}.
First, taking $G_1=\mathbb{R}^2$ and $G_2=SL_2(\mathbb{R})$, one
gets a matched pair of l.c. groups. Then one can show that the
corresponding split extension is a Kac algebra, which is
non-amenable, since $G_2$ is known to be a non-amenable l.c.
group.

Second, taking $G_1=\{(x,z)\vert x\in \mathbb{R},x\neq
0,z\in\mathbb{C}\}$ with multiplication $(x,z)(y,u)=(xy,z+xu)$,
and $G_2$ a double cover of $SU(1,1)$, one gets again a matched
pair of l.c. groups. Then one can show that the corresponding
split extension is not a Kac algebra. This l.c. quantum group is
non-compact, non-discrete and non-amenable, since $G_2$ is known
to be a non-amenable l.c. group.
\subsection{Kac exact sequence \cite{BSV1}}
We have seen above that one can associate an abelian group of
extensions with any matched pair $(G_1,G_2\subset G)$ of l.c.
groups. In the case of finite groups, Kac constructed in \cite{K1}
an {\it exact sequence} which allows to calculate the above group
of extensions in terms of usual cohomology groups of $G_1$, $G_2$
and $G$ with coefficients in the trivial module
$\mathbb{T}=\{z\in\mathbb{C}\vert \vert z\vert=1\}$. A similar
exact sequence was constructed in \cite{BSV1} for cohomology of
l.c. groups with coefficients in any Polish $G$-module $A$, but
here we discuss very briefly only the case $A=\mathbb{T}$.

First of all, let us precise, which kind of group cohomology we
deal with. Let us consider the cochain complex
$(L(\Gamma^n,\mathbb{T}))_n$, where $n=0,1,2,...,\
\Gamma=G_1,G_2,G,\ \Gamma^n$ is the Cartesian product of $n$
copies of $\Gamma,\ G^0$ is just a single point,
$L(\Gamma^n,\mathbb{T})$ is the set of (equivalence classes of)
Borel functions from $\Gamma^n$ to $\mathbb{T}$. The coboundary
operator $d : L(\Gamma^n,\mathbb{T})\to
L(\Gamma^{n+1},\mathbb{T})$ will be defined as follows. Let us
write face operators $\partial_i : \Gamma^{n+1}\to \Gamma^n$ :
$$
\partial_i(g_0,...,g_n)=\begin{cases} (g_1,...,g_n)\ \text{if}\
i=0,\cr (g_0,...,g_{i-1}g_i,...,g_n)\ \text{if}\ i=1,...,n,\cr
(g_0,...,g_{n-1})\ \text{if}\ i=n+1.\end{cases}
$$
Then, consider $ d_i : L(\Gamma^n,\mathbb{T})\to
L(\Gamma^{n+1},\mathbb{T})$ :
$$
(d_i F)(\bar{g})= \begin{cases} g_0\cdot F(\partial_0\bar{g}),\
\text{if}\ i=0,\cr F(\partial_i\bar{g}),\ \text{if}\
i=1,...,n+1,\end{cases}
$$
where $\overline{g}=(g_0,...,g_n)$. Finally,
$$
d=\Sigma_{i=0}^{n+1}(-1)^{-1}d_i.
$$
By definition, the {\it measurable cohomology} of the l.c. group
$\Gamma$ with coefficients in $\mathbb{T}$ is the cohomology of
the above cochain complex.

{\it Kac cohomology} $H(m.p.,\ \mathbb{T})$ of the matched pair
$(G_1,G_2\subset G)$ is defined as cohomology of certain complex
whose "building blocks" are $L(\Gamma_{p,q},\mathbb{T})$, the sets
of (equivalence classes of) Borel functions from $\Gamma_{p,q}$ to
$\mathbb{T}$, where $\Gamma_{p,q}$ is a closed subspace of
$G_1^{p(q+1)}\times G_2^{(p+1)q}$ defined in an inductive way
starting from $\Gamma_{p,0}=G_1^p,\ \Gamma_{0,q}=G_2^q,\
p,q=0,1,2,...$; typical example: $\Gamma_{1,1}=\{(g,h,s,t)\in
G_1\times G_1\times G_2\times G_2\vert sg=ht\}$. The coboundary
operator of this complex is constructed in terms of the above $d$
(see \cite{BSV1}). In particular, the group of extensions of the
matched pair $(G_1,G_2\subset G)$ is precisely the Kac
2-cohomology group $H^2(m.p.,\ \mathbb{T})$.

\cite{BSV1}, Corollary 4.5 claims that this Kac cohomology
$H(m.p.,\ \mathbb{T})$ satisfies the following long exact
sequence:
$$
0\to\mathbb{T}\to\mathbb{T}\oplus\mathbb{T}\to H^0(m.p.,\
\mathbb{T})\to H^1(G,\ \mathbb{T})\to H^1(G_1,\ \mathbb{T})\oplus
H^1(G_2,\ \mathbb{T})\to
$$
$$
\to H^1(m.p.,\ \mathbb{T}) \to H^2(G,\ \mathbb{T})\to H^2(G_1,\
\mathbb{T})\oplus H^2(G_2,\ \mathbb{T})\to H^2(m.p.,\
\mathbb{T})\to
$$
$$
\to H^3(G,\ \mathbb{T})\to H^3(G_1,\ \mathbb{T})\oplus H^3(G_2,\
\mathbb{T})\to ...
$$
Here $H^k(\Gamma,\mathbb{T}),\ k=1,2,...$ are the above mentioned
measurable cohomology groups of $\Gamma=G_1,G_2,G$ with
coefficients in $\mathbb{T}$.

This {\it Kac exact sequence} shows that in order to calculate the
group of extensions of the matched pair $(G_1,G_2\subset G)$, it
suffices to calculate $H^n(\Gamma,\mathbb{T})$ for $n=2,3$ and
$\Gamma=G_1,G_2,G$. This can be done for certain class of l.c.
groups, using the techniques proposed by D. Wigner. In particular,
in \cite{VV2} we computed the groups of extensions of matched
pairs of low dimensional Lie groups, passing to the corresponding
matched pairs of Lie algebras, and then performing  the
exponentiation of the results obtained on the Lie algebra level,
as it was explained in Section 3. But this last operation is quite
non-trivial, and it was precisely justified only in \cite{BSV1}.

\begin{remark} a) A cohomology theory for extensions of Lie
algebras was developed by A. Masuoka in \cite{Mas2}.

b) An example of calculation of the group of extensions for a
matched pair of l.c. groups via the Kac measurable cohomology can
be found in \cite{CKP2}.
\end{remark}
\end{section}


\begin{thebibliography}{A}

\bibitem{AW} H. Araki \& J. Woods : A classification of factors, {\it Publ. RIMS Kyoto
Univ., Ser. A}, {\bf 4} (1968), 51-130.

\bibitem{B} S. Baaj : Repr{\'e}sentation r{\'e}guli{\`e}re du groupe quantique
des d{\'e}placements de Woronowicz,
{\it Ast{\'e}risque}, {\bf 232} (1995), 11-48.

\bibitem{BS} S. Baaj \& G. Skandalis : Unitaires
multiplicatifs et dualit{\'e} pour les produits  crois{\'e}s de
${C}^*$-alg{\`e}bres, {\it Ann. Sci. ENS}, {\bf 26} (1993),
425-488.

\bibitem{BS1} S. Baaj \& G. Skandalis : Transformations pentagonales,
{\it C.R. Acad. Sci. Paris, S{\'e}r I}, {\bf 327} (1998), 623-628.

\bibitem{BSV} S. Baaj, G. Skandalis \& S. Vaes : Non-semi-regular quantum groups
coming from number theory, {\it Commun. Math. Phys.}, {\bf 235}
(1) (2003), 139-167.

\bibitem{BSV1} S. Baaj, Skandalis \& S. Vaes : Measurable Kac cohomology for
bicrossed products, {\it Trans. Amer. Math. Soc.}, {\bf 357}
(2005), 1497-1524.

\bibitem{BV2} S. Baaj \& S. Vaes : Double crossed products of locally compact
quantum groups, {\it Journal de l'Institut de Math. de Jussieu},
{\bf 4} (2005), 135-173.

\bibitem{BC} J.-B. Bost \& A. Connes : Hecke algebras, type III factors and
phase transitions with spontaneous symmetry breaking in number
theory, {\it Selecta Math. (N.S.)}, {\bf 1} (1995), 411-457.

\bibitem{BZ} F. Boca \& A. Zaharescu : Factors of type III and distribution of
prime numbers, {\it Proc. London. Math. Soc.}, {\bf 80} (3)
(2000), 145-178.

\bibitem{C} A. Connes : Une classification des facteurs de type III,
{\it Annales sci. de l'ENS, 4$^{e}$ s\'erie}, {\bf 6}, n.2 (1973),
133-252.

\bibitem{CKP1} Yu. A. Chapovsky, A. A. Kalyuzhnyi \& G. B. Podkolzin :
On 2+2 Locally Compact Quantum Groups, {\it Proc. of the Institute
of Mathematics of NAS of Ukraine}, {\bf 250} (3) (2004),
1064-1070.

\bibitem{CKP2} Yu. A. Chapovsky, A. A. Kalyuzhnyi \& G. B. Podkolzin :
On the group extensions for the bicrossed product construction for
a locally compact group, to appear.

\bibitem{DQV} P. Desmedt, J. Quaegebeur \& S. Vaes : Amenability
and the bicrossed product construction, {\it Illinois Journal of
Math.}, {\bf 46} (2002), 1259-1277.

\bibitem{D} V.G. Drinfeld : Quantum groups, In {\it Proceedings of ICM}, (Berkeley, Calif.,
1986), Amer. Math. Soc., Providence, RI, 1987, 798-820.

\bibitem{ES} M. Enock \& J.-M. Schwartz : Kac algebras and duality of locally compact
groups,  Springer, Berlin, 1992.

\bibitem{ES1} M. Enock \& J.-M. Schwartz : Alg\`ebres de Kac moyennables, {\it Pacific J. Math.},
{\bf 125} (2) (1986), 363-379.


\bibitem{F} P. Fima : Locally compact quantum groups whose algebras are factors (in preparation).

\bibitem{HR}  E.  Hewitt \&  K.A.  Ross  :  Abstract  Harmonic Analysis, v.II,
Springer-Verlag, Berlin, 1970.

\bibitem{JS}  V. Jones \& V. S. Sunder :  {Introduction to
Subfactors}, London Math. Society Lecture Notes Series, {\bf 234},
University Press, Cambridge, 1997.

\bibitem{K} G.I. Kac : Generalization of the group principle of duality,
{\it Soviet Math. Dokl.}, {\bf 2} (1961), 581-584. Ring-groups and
the principle of duality, I, II, {\it Trans. Moscow Math. Soc.},
{\bf 12} (1963), 291-339; {\bf 13} (1965), 94-126.

\bibitem{K1} G.I. Kac : Extensions of groups to ring-groups, {\it Math. USSR
Sbornik}, {\bf 5} (1968), 451-474.

\bibitem{KP1} G.I. Kac \& V.G. Paljutkin : Finite ring groups, {\it Trans.
Moscow Math. Soc.}, {\bf 5} (1966), 251-294.

\bibitem{KP2} G.I. Kac \& V.G. Paljutkin : An example of a ring group of order
eight (Russian), {\it Soviet Math. Surveys}, {\bf 20, n.5} (1965),
268-269.

\bibitem{KP3} G.I. Kac \& V.G. Paljutkin : Example of a ring group
generated by Lie groups (Russian), {\it Ukr. Math. J.}, {\bf 16,
n.1} (1965), 268-269.

\bibitem{KK} E. Koelink \& J. Kustermans : A locally compact quantum group
analogue of the normalizer of $SU(1,1)$ in $SL(2,C)$, {\it Comm.
Math. Phys.}, {\bf 233} (2003) 231-296.

\bibitem{KK1} E. Koelink \& J. Kustermans : Quantum $\widetilde{SU(1,1)}$
and its Pontryagin dual, In {\it Locally Compact Quantum Groups
and Groupoids. Proceedings of the Meeting of Theoretical
Physicists and Mathematicians, Strasbourg, February 21-23, 2002},
Ed. L. Vainerman, IRMA Lectures on Mathematics and Mathematical
Physics, {\bf 2}, Walter de Gruyter, Berlin, New York  (1) (2003),
49-77.

\bibitem{KV1} J. Kustermans \& S. Vaes : Locally compact
quantum groups. {\it Ann. Scient. Ec. Norm. Sup.} {\bf 33} (6)
(2000), 837--934.

\bibitem{KV2} J. Kustermans \& S. Vaes : Locally compact quantum groups in
the von Neumann algebraic setting, {\it Math. Scand.}, {\bf 92}
(1) (2003), 68-92.

\bibitem{KV3} J. Kustermans \& S. Vaes : A simple
definition  for locally compact quantum groups. {\it C.R. Acad.
Sci., Paris, S{\'e}r. I} {\bf 328} (10) (1999), 871--876.

\bibitem{KVVVW} J. Kustermans, S. Vaes, L. Vainerman, A. Van
Daele \& S.L. Woronowicz : Locally Compact Quantum Groups. In {\it
Lecture Notes for the school on Noncommutative Geometry and
Quantum Groups in Warsaw (17-29 September 2001)}. {\it Banach
Centre Publications}, to appear.


\bibitem{Mac} G.W. Mackey : {Products of subgroups and projective
multipliers}, {\it Colloquia Mathematica Societatis J\`anos
Bolyai}, {\bf 5}. Hilbert space operators. Tihany (Hungary)
(1970), 401-413.

\bibitem{Maj1} S. Majid : Physics for algebraists: Non-commutative and
non-cocommutative Hopf algebras by a bicrossproduct construction.
{\it J. Algebra} {\bf 130} (1990), 17--64.

\bibitem{Maj2} S. Majid : More examples of bicrossproduct and double
cross product Hopf algebras. {\it Isr. J. Math.} {\bf 72} (1990),
133--148.

\bibitem{Maj3} S. Majid : Hopf-von Neumann algebra bicrossproducts, Kac
algebra bicrossproducts, and the classical Yang-Baxter equations.
{\it J. Funct. Anal.} {\bf 95} (1991), 291--319.

\bibitem{Majbook} S. Majid : Foundations of Quantum Group Theory,
Cambridge University Press, 1995.

\bibitem{MN} T. Masuda \& Y. Nakagami : A von Neumann Algebra
Framework for the Duality of the Quantum Groups. {\it Publ. RIMS,
Kyoto Univ.}, {\bf 30}, n.5 (1994), 799-850.

\bibitem{MNW} T. Masuda, Y. Nakagami \& S.L. Woronowicz : A
$C^*$-Algebraic Framework for Quantum Groups. {\it Intl. Journal
of Math.}, {\bf 14}, n.9 (2003), 903-1001.


\bibitem{Mas2} A. Masuoka :  Extensions of Hopf Algebras and Lie
Bialgebras. {\it Trans. of the AMS} {\bf 352}, No. 8 (2000),
3837--3879.




\bibitem{PW} P. Podle\'s \& S.L. Woronowicz : {Quantum deformation of
Lorenz group}, {\it Comm. Math. Phys.}, {\bf 130} (1990), 381-431.

\bibitem{Ruan} Zh.-J. Ruan : {Amenability of Hopf von Neumann algebras and Kac
algebras}, {\it J. Funct. Anal.}, {\bf 139, n.2} (1996), 466-499.

\bibitem{Sauv} J-L. Sauvageot : Sur le type du produit crois\'e d'une
alg\`ebre de von Neumann par un groupe localement compact, {\it
Bulletin de la Soc. Math. France}, {\bf 105} (1977), 349-346.

\bibitem{Skand} G. Skandalis : Duality for locally compact 'quantum groups'
(joint work with S.~Baaj), Mathematisches Forschungsinstitut
Oberwolfach, Tagungsbericht 46/1991, C*-algebren, 20.10 --
26.10.1991, p.~20.

\bibitem{Stra} S. Stratila : Modular Theory in Operator Algebras.
{\it Abacus Press, Tunbridge Wells, England} (1981).

\bibitem{Su} C. Sutherland : The analysis of the regular representation of a
non-unimodular group, {\it Pacif. Journal of Math.}, {\bf 79} (1)
(1978), 225-250.

\bibitem{Tak} M. Takeuchi : Matched pairs of groups and bismash products of Hopf
algebras, {\it Comm. Algebra}, {\bf 9} (1981), 841-882.

\bibitem{Tom} R. Tomatsu : {Amenable discrete quantum groups},
Preprint {\it arXiv:math.QA/0302222}, (2003), 22p.

\bibitem{SV2} S. Vaes : Examples of locally compact quantum groups through the
bicrossed product construction. {\it Proceedings of the XIIIth
Int. Conf. Math. Phys. London, 2000}, Eds. A. Grigoryan, A. Fokas,
T. Kibble \& B. Zegarlinski, Intl. Press of Boston, Sommerville MA
(2001), 341-348.

\bibitem{VV1} S. Vaes \& L. Vainerman : Extensions of locally compact quantum groups
and the bicrossed product construction, {\it Adv. in Math.}, {\bf
175} (1) (2003), 1-101.

\bibitem{VV2} S. Vaes \& L. Vainerman : On low dimensional locally compact quantum
groups, In {\it Locally Compact Quantum Groups and Groupoids.
Proceedings of the Meeting of Theoretical Physicists and
Mathematicians, Strasbourg, February 21-23, 2002}, Ed. L.
Vainerman, IRMA Lectures on Mathematics and Mathematical Physics,
{\bf 2}, Walter de Gruyter, Berlin, New York  (1) (2003), 127-187.


\bibitem{VD2} A. Van Daele : Multiplier Hopf Algebras,
{\it Trans. Amer. Math. Soc.} {\bf 342} (1994), 917--932.

\bibitem{VD3} A. Van Daele : An algebraic framework for group
duality, {\it Adv. in Math.} {\bf 140} (1998), 323--366.

\bibitem{VW} A. Van Daele \& S.L. Woronowicz : Duality for the quantum $E(2)$ group,
{\it Pac. J. Math.}, {\bf 173, n.2} (1996), 375-385.

\bibitem{W1} S.L. Woronowicz : Compact matrix pseudogroups, {\it Comm. Math. Phys.}
{\bf 111} (1987), 613-665.

\bibitem{W2} S.L. Woronowicz : Twisted $SU(2)$ group. An example of a non-commutative
differential calculus, {\it Publ. RIMS, Kyoto University}, {\bf
23} (1987), 117-181.

\bibitem{W3} S.L. Woronowicz : Tannaka-Krein duality for compact matrix
pseudogroups, Twisted $SU(N)$ groups. {\it Invent. Math.} {\bf 93}
(1988), 35-76.

\bibitem{W4} S.L. Woronowicz : Quantum $E(2)$ group and its Pontryagin dual, {\it Lett.
Math. Phys.} {\bf 23} (1991), 251-263.

\bibitem{W5} S.L. Woronowicz : From multiplicative unitaries to quantum
groups, {\it Int. J. Math.}, {\bf 7}, n.1 (1996), 127-149.

\bibitem{W6} S.L. Woronowicz : Quantum $"az+b"$
group on complex plane, {\it Int. J. Math.}, {\bf 12}, n.4 (2001),
461-503.

\bibitem{W7} S.L. Woronowicz \& S. Zakrzevski : Quantum $"ax+b"$
group, {\it Review in Math. Phys.}, {\bf 14}, n.7 (2002), 797-828.

\bibitem{W8} S.L. Woronowicz : Extended $SU(1,1)$ quantum group. Hilbert space
level, {\it Preprint KMMF} (in preparation).
\end{thebibliography}
\end{document}